\documentclass[12pt]{article}
\usepackage{amsmath}
\usepackage{amsthm}
\usepackage{amssymb}
\usepackage{bigints}
\usepackage{graphicx}
\usepackage{booktabs}
\usepackage{threeparttable}
\title{On the speed of the one-dimensional polymer in the large range regime}
\author{Chien-Hao Huang\footnote{chhuang@ncts.ntu.edu.tw, National Center for Theoretical Sciences, No. 1 Sec. 4 Roosevelt Rd.,
National Taiwan University
Taipei, 10617, Taiwan} \\
NCTS}
\date{ }

\newtheorem{theorem}{Theorem}[section]
\newtheorem{lemma}[theorem]{Lemma}

\newtheorem{corollary}[theorem]{Corollary}

\numberwithin{equation}{section}

\newcommand{\SortNoop}[2]{#2}

\begin{document}
\maketitle
\fontsize{12}{18pt}\selectfont

\begin{abstract}
    We consider a Hamiltonian involving the range of the simple random walk and the Wiener sausage so that the walk tends to stretch itself. This Hamiltonian can be easily extended to the multidimensional cases, since the Wiener sausage is well-defined in any dimension. In dimension one, we give a formula for the speed and the spread of the endpoint of the polymer path. Also, we provide the CLT. It can be easily showed that if the self-repelling strength is stronger, the end point is going away faster. This strict monotonicity of speed has not been proven in the literature for the one-dimensional case.     
\end{abstract}


\section{Introduction}

\subsection{The model and main results}

A polymer consists of monomers. Monomers have tendency to repel each other because two monomers can not occupy the same site. This phenomenon is called the excluded-volume-effect. There is a probabilistic way to model this physical phenomenon (cf. Madras and Slade \cite{MS93} section 2.2). We use the $d$-dimensional simple random walk $\{S_n\})_{\mathbb{N}\cup \:0}$ under the probability measure $P$ to represent the position of monomers and $S_n$ to represent the end-point of the polymer chain with length $n$. 
$S_0=0$ and $S_n=\sum_{i=1}^n X_i$, where $(X_i)_{i\in\mathbb{N}}$ is a sequence of independent and identically distributed (i.i.d.) random variables. The distribution of $X_i$'s is
\[P(X_1=x)=\left\{\begin{array}{l} \frac{1}{2d}, \; x\in \mathbb{Z}^d \; \mbox{with} \; ||x||=1,\\
                              0,  \;\; \mbox{otherwise}.
\end{array}\right.\]
The random process $(S_n)_{n\in\mathbb{N}\cup 0}$ is called the \textit{simple random walk} (SRW) on $\mathbb{Z}^d$.
Suppose that the end-point has scale $\alpha_n$, the local density of monomers will be $\frac{n}{\alpha_n^d}$. The self-repelling energy is approximately
\begin{equation}\label{heuristic}
\exp\left(Energy\right) \approx \exp\left(-\sum_{x\in \mathbb{Z}^d}\left(\frac{n}{\alpha_n^d}\right)^2\mathbf{1}_{x \text{ is occupied}}\right) \approx \exp\left( -\alpha_n^d \times \left(\frac{n}{\alpha_n^d}\right)^2 \right).
\end{equation}
On the other hand, by the local limit theorem of the simple random walk,
\begin{equation}
P(|S_n| = \alpha_n)\approx \exp(-C\alpha_n^2/n).
\end{equation}
Let
\begin{equation}
\frac{n^2}{\alpha_n^d} = \frac{\alpha_n^2}{n},
\end{equation}
we get $\alpha_n= n^{\frac{3}{d+2}}$. It is expected that $|S_n|\sim n^{\frac{3}{d+2}}$ for $d=1,2,3$, and $|S_n|\sim n^{1/2}$ for $d\geq 4$ with a logarithmic correction when $d=4$ under the self-repelling phenomenon.\\

In this paper, we propose the following Hamiltonian
\begin{equation} G_n := \frac{n^2}{R_n},\end{equation}
where $R_n$ is the number of sites occupied by the walk up to time $n-1$, that is,
\begin{equation}
R_n:= \# \{x: \exists\: i,\; S_i=x, \;0\leq i\leq n-1  \}.
\end{equation}
Fix $n\in\mathbb{N}$ and a parameter $\beta\in (0,\infty)$, denote
\begin{equation}
Z_n^{G} := E \left(\exp \left(-\beta G_n \right)\right)
\end{equation}
and
\begin{equation}
Z_n^{G}(A) := E \left(\mathbf{1}_A\exp \left(-\beta G_n \right)\right).
\end{equation}
\noindent
The polymer measure is then defined by
\begin{equation} P_n^{G} (S):= \frac{e^{-\beta G_n(S)}}{Z_n^{G}} P(S). \end{equation}
\noindent
$\beta$ is called the strength of the self-repellence. This polymer measure favors the event ``the polymer has large range''.

Let $I_d(x):=\lim_{n\rightarrow \infty}\frac{-1}{n}\log P\{R_n\geq xn\}$ and $I(x):=I_1(x)=\frac{1}{2}(1+x)\log(1+x)+\frac{1}{2}(1-x)\log(1-x)$. The following are our main results for the one-dimensional discrete setting.
\begin{theorem}  
(i) For $\beta>0$,
\begin{equation}
\lim_{n\rightarrow \infty}\frac{1}{n} \log Z_n^{G}= g^*(\beta),
\end{equation}
\noindent
where
\begin{equation}\label{inf}
g^*(\beta):=-\inf_{c\in [\tilde{c}(\beta) , 1]} \left\{ \frac{\beta}{c} + I(c) \right\}
\end{equation}
and $\tilde{c}(\beta)=\frac{\beta}{\beta+ \log 2}$. \\
(ii) The infimum of \eqref{inf} is obtained at $c^*(\beta)$, 
where $c^* (\beta)$ is the solution of 
\begin{equation}
\beta =c^2 I'(c)=\frac{c^2}{2} \log \left( \frac{1+c}{1-c}\right).
\end{equation}
Note that $c^*$ is strictly monotone, $\beta^{-1/3} c^*(\beta)\rightarrow 1$ as $\beta \rightarrow 0$ and $e^{2\beta}(1-c^*(\beta))\rightarrow 2$ as $\beta \rightarrow \infty$.
$g^*(\beta)$ can be written as
\begin{equation}
g^*(\beta) = -c^*(\beta) \log\frac{1+c^*(\beta)}{1-c^*(\beta)} -\frac{1}{2}\log(1-c^*(\beta)^2).
\end{equation}
Furthermore, $\beta^{-2/3} g^*(\beta)\rightarrow -\frac{3}{2}$ as $\beta \rightarrow 0$ and $g^*(\beta)+\beta \rightarrow -\log 2$ as $\beta \rightarrow \infty$.
\end{theorem}

\begin{theorem}(LLN and LDP)  
For $\beta>0$, define
\begin{equation}
P_n^{G,+}(\cdot)= P_n^G(S_n/n\in \cdot| S_n >0)
\end{equation}
Then $(P_n^{G,+})_{n\in \mathbb{N}}$ satisfies a lage deviation principle (LDP) on $[0,1]$ with rate $n$ and with rate function $I^{\beta}(\theta)$ 
\[-I^{\beta}(\theta) =\left\{\begin{array}{ll}
 -\frac{\beta}{\theta} -I(\theta) -g^{*}(\beta), &  c^*(\frac{\beta}{2})\leq \theta,  \\
 -\frac{\beta}{\tilde{r}}-I(2\tilde{r}-\theta) -g^{*}(\beta), &   \theta  < c^*(\frac{\beta}{2}),
\end{array}\right. \]
where $\tilde{r}= \tilde{r}_{\beta}(\theta)$ is the positive solution of $\beta =2 r^2 I'(2r-\theta)$. Moreover, $I^{\beta}(\theta)$ has the unique $0$ at $c^*(\beta)$
\end{theorem}

\begin{theorem}(CLT) 
$\forall \; C \in \mathbb{R}$,
\begin{equation}
\lim_{n\rightarrow \infty}  P^G_n\left( \frac{ S_n-c^*(\beta)n}{\sigma^*(\beta)\sqrt{n}}\leq C | S_n > 0\right) =\Phi(C),
\end{equation}
where $\frac{1}{ \sigma^{*2} (\beta)} = \left( \frac{\beta}{\theta}+I(\theta)\right)''|_{\theta=c^*(\beta)}= \frac{2\beta}{c^{*3}(\beta)} + \frac{1}{1-c^*(\beta)^2}$.
$\sigma^*(\beta)\rightarrow \frac{1}{\sqrt{3}}$ as $\beta \rightarrow 0$ and $e^{\beta}\sigma^*(\beta)\rightarrow 2$ as $\beta \rightarrow \infty$. $\Phi$ is the cumulative distribution function of the standard normal random variable.
\end{theorem} 

Our Hamiltonian can be extended to the continuous setting easily. Let $B_t$ be the $d-$dimensional Brownian motion, the analogous Hamiltonian is 
\begin{equation}G_t:=\frac{t^2}{R^1_t},\end{equation}  
where $R^1_t$ is the Wiener sausage with radius 1 in $d$-dimension. Note that
\begin{equation}
R_t^{\epsilon}:= Volume\left(\bigcup_{s\leq t} Ball_{\epsilon}(B_s) \right).
\end{equation}
Denote
\begin{equation}
Z_t^{G} := E \left(\exp \left(-\beta G_t \right)\right)
\end{equation}
and
\begin{equation}
Z_t^{G}(A) := E \left(\mathbf{1}_A\exp \left(-\beta G_t \right)\right).
\end{equation}
\noindent
The polymer measure is defined by
\begin{equation} dP_t^{G} := \frac{e^{-\beta G_t}}{Z_t^{G}} dP. \end{equation}
$R_t^1$ is finite for any $t$ and in any dimension. Thus, the polymer measure is well-defined in any dimension.

Let $J_d(x):=\lim_{n\rightarrow \infty}\frac{-1}{t}\log P\{R_t^1\geq xt\}$ and $J(x):=J_1(x)=\frac{x^2}{2}$.
We have our main results for the one-dimensional continuous setting as following.

\begin{theorem}  
(i) For $\beta>0$,
\begin{equation}
\lim_{t\rightarrow \infty}\frac{1}{t} \log Z_t^{G}= g^{**}(\beta),
\end{equation}
\noindent
where
\begin{equation}\label{inf2}
g^{**}(\beta):=-\inf_{c\in [\tilde{\beta_1},\infty)} \left\{ \frac{\beta}{c} + J(c) \right\}
\end{equation}
and $\tilde{\beta_1} :=\frac{1}{2}\beta^{1/3}$.
\noindent
(ii) The infimum of \eqref{inf2} is obtained at $c^{**} (\beta) = \beta^{1/3}$ and $g^{**} (\beta) = -\frac{3}{2}\beta^{2/3}$. Moreover, $Z^G_t \sim \frac{8}{\sqrt{3}} e^{g^{**} (\beta)t}$.
\end{theorem}

\begin{theorem}(LLN and LDP)  
For $\beta>0$, define
\begin{equation}
P_t^{G,+}(\cdot)= P_t^G(B_t/t\in \cdot| B_t >0)
\end{equation}
Then $(P_t^{G,+})_{t>0}$ satisfies a lage deviation principle (LDP) on $[0,1]$ with rate $n$ and with rate function $J^{\beta}(\theta)$ 

\[-J^{\beta}(\theta) =\left\{\begin{array}{ll}
 -\frac{\beta}{\theta} - J(\theta) -g^{**}(\beta), &  \sqrt[3]{\frac{\beta}{2}}\leq \theta,  \\
 -\frac{\beta}{\bar{r}}-J(2\bar{r}-\theta) -g^{**}(\beta), &   \theta  < \sqrt[3]{\frac{\beta}{2}},
\end{array}\right. \]
where $\bar{r}= \bar{r}_{\beta}(\theta)$ is the positive solution of $\beta =2r^2 J'(2r-\theta)$. $(\sqrt[3]{\frac{1}{2}} \approx 0.79)$
\end{theorem}

\begin{theorem}(CLT) 
$\forall \; C \in \mathbb{R}$,
\begin{equation}
\lim_{t\rightarrow \infty}  P^G_t\left( \frac{ B_t-c^{**}(\beta)t}{\sigma^{**}(\beta)\sqrt{t}}\leq C | B_t> 0\right) =\Phi(C),
\end{equation}
where $\sigma^{**} (\beta) =\frac{1}{\sqrt{3}}$.
\end{theorem}

Remark. In higher dimensions $d\geq 2$, Theorem 1.1 (i) and 1.4 (i) have analogous results since $I_d$ and $J_d$ are well studied in Hamana and Kesten \cite{HK01}.
\begin{corollary}  
For $\beta>0$,
\begin{equation}
\lim_{n\rightarrow \infty}\frac{1}{n} \log Z_n^{G}= g^*(\beta),
\end{equation}
\noindent
where
\begin{equation}
g^*(\beta):=-\inf_{c\in [\tilde{c}_d(\beta) , 1]} \left\{ \frac{\beta}{c} + I_d(c) \right\}
\end{equation}
and $\tilde{c}_d(\beta)=\frac{\beta}{\beta+ \log 2d}$.

For the continuous setting,
\begin{equation}
\lim_{t\rightarrow \infty}\frac{1}{t} \log Z_t^{G}= g^{**}(\beta),
\end{equation}
\noindent
where
\begin{equation}
g^{**}(\beta):=-\inf_{c\in [\tilde{\beta_d},\infty)} \left\{ \frac{\beta}{c} + J_d(c) \right\}
\end{equation}
and $\tilde{\beta_d} :=\frac{1}{2}\left(\frac{\beta}{w_{d-1}}\right)^{1/3}$. $w_{d-1}$ is the volume of the unit ball in $(d-1)$-dimension. Set also $w_0=1$.
\end{corollary}

{\bf Remark.} When $d=2$, Corollary 1.7 implies that $\liminf R_n/n$ is at least $\frac{\beta}{\beta +\log 4}$ under the polymer measure. In another hand, $R_n^{(\alpha)}/n$ of a two-dimensional $\alpha$-stable random walk $S^{(\alpha)}_n$ has a positive limit if and only if $\alpha \leq 4/3$ $\cite{LR91}$. Recall that $S_n^{(\alpha)}$ is of order $n^{1/\alpha}$. We expect that the endpoint behaviour of the two-dimensional self-repellent polymer is at least of order $n^{3/4}$.

\subsection{Literature and discussions}

In the literature, the Domb-Joyce model, relying on the random walk setting, and the Edwards model, which is based on Brownian motion, both discuss the self-repellence concerning the number of self-intersections.
The Hamiltonian
\begin{equation}
H_n(S) :=\sum_{i,j=0; \; i\neq j}^{n-1} \textbf{1}_{S_i=S_j}= \sum_{x\in\mathbb{Z}^d} \ell_n^2(x) -n
\end{equation}
is the self-intersection local time up to time $n-1$, and
$$\ell_n(x)=\#\{0\leq i\leq n-1 : S_i=x\}, \;\;\; x\in\mathbb{Z}^d,$$
is the local time at site $x$ up to time $n-1$. 
The polymer measure is then defined by
\begin{equation} P_n^{H} (S):= \frac{e^{-\beta H_n(S)}}{Z_n^{H}} P(S), \end{equation}
where
\begin{equation}
Z_n^{H} := E \left(\exp \left(-\beta H_n \right)\right).
\end{equation}
The path receives a penalty $e^{2\beta}$ when the path self-intersects itself. This model is also called the \textit{weakly self-avoiding walk} (WSAW). For more details about this model, one can see the monograph by den Hollander \cite{dH09}. Readers can see that $\beta=0$, $P^H_n$ is the measure for SRW and $\beta=\infty$, $P^H_n$ is the uniform measure for the \textit{self-avoiding walk} (SAW), namely, no self-intersection is allowed. The discussion for SAW has been lasting for a long time, please see the monograph \cite{MS93} for the earlier work and the recent survey Bauerschmidt et al. \cite{BDGS12}. People believe that the WSAW $(0<\beta<\infty)$ and SAW ($\beta=\infty$) are in the same universality, i.e., the exponent for both model are the same: $\frac{3}{d+2}\wedge \frac{1}{2}$. 
Define
\begin{equation}
\hat{H}_n := \sum_{x\in \mathbb{Z}} \ell_n^2(x).
\end{equation}
Since the polymer measure will not be changed by taking away constants in the Hamiltonian $H_n$, the analysis of $\hat{H}_n$ is the same as $H_n$. Our Hamiltonian $G_n$ is weaker than $\hat{H}_n$. For ``weaker'' we mean that
\begin{equation}
\left[ \:\sum_{x\in\mathbb{Z}^d} \ell_n^2(x) \right] \cdot \left[\: \sum_{x\in\mathbb{Z}^d} \mathbf{1}_{\ell_n(x)>0} \right] \geq \left[\: \sum_{x\in\mathbb{Z}^d} \ell_n(x)\mathbf{1}_{\ell_n(x)>0} \right]^2 
= \left[\: \sum_{x\in\mathbb{Z}^d} \ell_n(x) \right]^2
=n^2.
\end{equation}
We have
\begin{equation}
\hat{H}_n \geq \frac{n^2}{R_n}=G_n.
\end{equation}
Concerning $G_n -n= n\left(\frac{n-R_n}{R_n}\right)$. When $\beta=0$, $P^G_n$ is SRW. On the other hand, $\beta=\infty$, the event $\{R_n=n\}$ survives under $P^G_n$, but $\{R_n=k\}$ for $1\leq k< n$ all die out.  
$G_n$ should be in the same universality as WSAW and SAW.

For the continuous setting,
\begin{equation}
H_t(B):= \int_{0}^t  ds \int_0^t du \;\;\delta_{0}(B_s -B_u),
\end{equation}
$H_t$ is infinity when the dimension is higher than one. Past results used truncations to obtain the polymer measure as a weak limit. See van der Hofstad \cite{vdH98} and the reference therein. However, our Hamiltonian $G_t$ is well-defined for any $t$ in any dimension, so is the polymer measure $P^G_t$.

The results for one-dimensional WSAW are rich. Greven and den Hollander \cite{GdH93} proved the LDP for the Domb-Joyce model and K{\"o}nig \cite{Ko96} gave the CLT result. For the one-dimensional Edwards model, van der Hofstad et al. \cite{vdHdHK03} proved the LDP and van der Hofstad et al. \cite{vdHdHK97AOP} showed the CLT. In van der Hofstad \cite{vdH98JSP}, the author had numerical results for the speed and the variance of the speed in \cite{vdHdHK97AOP}, that is,
$c^{**}(\beta)\beta^{-1/3}\in [1.104,1.124]$ and $\sigma^{**}(\beta)\in [0.60,0.66]$, while we have $c^{**}(\beta)\beta^{-1/3}=1$ and $\sigma^{**}(\beta)=1/\sqrt{3}\fallingdotseq 0.577$ in Theorem 1.5 and 1.6.

In the rest of paper, Section 2.1 gives the proofs for Theorem 1.4-1.6 and Section 2.2 gives the proofs for Theorem 1.1-1.3. 

\section{Proofs} 

When $d=1$, we not only have the explicit formula of the density of $$R_t:=\max_{0\leq s\leq t} B_s - \min_{0\leq s\leq t} B_s = R^1_t-2,$$ but also have the joint density of $B_t$ and $R_t$. Therefore, we start the proofs for the continuous model. The following lemma plays an important role.
\begin{lemma}
Let $\{A_n\}$ and $\{B_n\}$ be two positive sequences, we have that $$\limsup_{n\rightarrow \infty} \frac{1}{n}\log (A_n+B_n) \leq \limsup_{n\rightarrow \infty} \frac{1}{n}\log A_n \vee \limsup_{n\rightarrow \infty} \frac{1}{n}\log B_n.$$
\end{lemma}
\noindent
$\mathbf{Proof}$ $\;$ Trivial. \qed

In the rest of paper, $f(t)\sim g(t)$ means that $f(t)/g(t)$ goes to $1$ as $t$ goes to $\infty$.

\subsection{The continuous model} 
The rate function for the large deviation principle for the Wiener sausage was discussed.
\begin{theorem} (Hamana and Kesten \cite{HK01}) 
\begin{equation}
J_d(x)=\lim_{n\rightarrow \infty}\frac{-1}{t}\log P\{R_t^1\geq xt\}
\end{equation}
exists in $[0,\infty)$ for all $x$. $J_d(x)$ is continuous on $[0,\infty)$ and strictly increasing on $[C_d,\infty)$, and for $d\geq 2$, $J_d(x)$ is convex on $[0,\infty)$. Furthermore,
\begin{equation}
  \begin{aligned}
J_d(x)=0 & \;\;\it{for}\;\; x\leq C_d,\\
0<J_d(x)<\infty & \;\;\it{for}\;\; C_d < x.
   \end{aligned}
\end{equation}
\end{theorem}
\noindent
For $d=1$, $J(x):=J_1(x)=\frac{x^2}{2}$ for $x\geq 0$. $C_d$ is the heat capacity of the unit ball for the $d$-dimension Brownian motion. The large deviation principle was also provided.\\

\noindent
$\mathbf{Proof \;\;of \;\;Theorem\;\; 1.4}$\\ 
In the case $d=1$, we have the explicit formula of the density of 
\begin{equation}
R_t:=\max_{0\leq s\leq t} B_s - \min_{0\leq s\leq t} B_s = R^1_t-2.
\end{equation} From Feller \cite{Fe51},
\begin{equation}\label{densityr}
P(R_t\in dr)= \frac{8}{\sqrt{ t}} \sum_{k=1}^{\infty} (-1)^{k-1} k^2 \phi\left(\frac{kr}{t^{1/2}}\right)dr,
\end{equation}
where $\phi(x)=\frac{1}{\sqrt{2\pi}} e^{-x^2/2}$. Notice that the series converges uniformly on $[a,\infty)$ for any $a>0$, but \textbf{not} on $[0,\infty)$. In the proof, we always consider the event that the variables are away from zero. 
We will find out soon that $R^1_t$ is of order $t$, thus, we use $R_t$ instead of $R^1_t$.
The argument of the proof is that we first match the order 
\begin{equation}
\frac{t^2}{r}=\frac{r^2}{t},
\end{equation}
so that we know the range is of order t. Later, we set $r=ct$ and do the change of variables. The first term in \eqref{densityr} is dominant. By Laplace's method, we get 
\begin{equation}
\frac{\beta}{c^2} =J'(c)=\left(\frac{c^2}{2}\right)^{'}=c \;\;\mbox{and}\;\; c=\beta^{1/3}. 
\end{equation}
Rigorously, we apply Theorem 2.2 and the Varadhan's lemma,
\begin{equation}
\lim_{t\rightarrow \infty}\frac{1}{t}\log E\left( e^{-\beta \frac{t^2}{R_t }} \mathbf{1}_{R_t\geq \frac{1}{2}\beta^{1/3} t} \right) = -\inf_{c\in [\frac{1}{2}\beta^{1/3} ,\infty) } \left\{ \frac{\beta}{c} + J(c) \right\} =-\frac{3}{2}\beta^{2/3}
\end{equation}
and
\begin{equation}
\limsup_{t\rightarrow \infty}\frac{1}{t}\log E\left( e^{-\beta \frac{t^2}{R_t}} \mathbf{1}_{R_t\leq \frac{1}{2}\beta^{1/3} t} \right) \leq
-2\beta^{2/3}. 
\end{equation}
By Lemma 2.1, (i) is proved. For $d\geq 2$, the volume of the sausage is larger than the volume of the sausage of the first coordinate. 
\begin{equation}
\liminf_{t\rightarrow \infty}\frac{1}{t}\log E\left( e^{-\beta \frac{t^2}{R_t^1 }} \right) \geq \lim_{t\rightarrow \infty}\frac{1}{t}\log E\left( e^{-\beta \frac{t^2}{w_{d-1}(R_t +2) }}  \right) = -\frac{3}{2} \left(\frac{\beta}{w_{d-1}}\right)^{2/3}
\end{equation}
and 
\begin{equation}
\limsup_{t\rightarrow \infty}\frac{1}{t}\log E\left( e^{-\beta \frac{t^2}{R_t^1}} \mathbf{1}_{R_t^1\leq \tilde{\beta}_d t} \right) \leq
-2 \left(\frac{\beta}{w_{d-1}}\right)^{2/3}.
\end{equation}
This proves the second part of Corollary 1.7. Moreover, 
\begin{equation}
\limsup_{t\rightarrow \infty}\frac{1}{t}\log E\left( e^{-\beta \frac{t^2}{R_t }} \mathbf{1}_{R_t\geq \frac{7}{4}\beta^{1/3} t} \right) \leq \limsup_{t\rightarrow \infty}\frac{1}{t}\log E\left( \mathbf{1}_{R_t\geq \frac{7}{4}\beta^{1/3} t} \right) 
=-\frac{49}{32}\beta^{2/3} .
\end{equation}
This means that we only need to consider $R_t/t$ between $\frac{1}{2}\beta^{1/3}$ and $\frac{7}{4}\beta^{1/3}$. Therefore, we only need to consider the first term in \eqref{densityr}. We also need these two bounds to compute the Taylor series of $-\left(\frac{\beta}{c} + \frac{c^2}{2}\right)$ at $\beta^{1/3}$ which is 
$$-\frac{3}{2}\beta^{2/3} -\frac{3}{2}(c-\beta^{1/3})^2 +\sum_{k=3}^{\infty} \frac{(-1)^{k+1}}{\beta^{(k+1)/3}} (c-\beta^{1/3})^k $$ for $|c- \beta^{1/3}| < \beta^{1/3}$. Now, we compute the asymptotic of $Z_t^G$.
\[\begin{array} {rcl}
Z_t^G & \sim &\frac{8\sqrt{t}}{\sqrt{2\pi}}  \int_{\frac{1}{2}\beta^{1/3}}^{\frac{7}{4}\beta^{1/3}}  e^{ - \frac{\beta}{c} t -\frac{c^2}{2}t } dc \\
	&=& \frac{8\sqrt{t}}{\sqrt{2\pi}} e^{ -\frac{3}{2}\beta^{2/3}t} \int_{\frac{1}{2}\beta^{1/3}}^{\frac{7}{4}\beta^{1/3}}  e^{  -\frac{3}{2}(c-\beta^{1/3})^2 t +\sum_{k=3}^{\infty} \frac{(-1)^{k+1}}{\beta^{(k+1)/3}} (c-\beta^{1/3})^k t} dc\\
    &=& \frac{8}{\sqrt{2\pi}} e^{ -\frac{3}{2}\beta^{2/3}t} \int_{-\frac{1}{2}\beta^{1/3}\sqrt{t}}^{\frac{3}{4}\beta^{1/3}\sqrt{t}}  e^{  -\frac{3}{2}c^2 +\sum_{k=3}^{\infty} \frac{(-1)^{k+1}}{\beta^{(k+1)/3}} \left(\frac{c}{\sqrt{t}}\right)^k t} dc\\    
	&=& \frac{8}{\sqrt{3}} e^{ -\frac{3}{2}\beta^{2/3}t}  \left(1+O(\frac{1}{t})\right).

\end{array}\]
One can see that those terms in the exponent with degree higher than 2 are of order $O\left(\frac{1}{t}\right)$, they are negligible. As a by-product, for the second order approximation of $R_t$ under $P^G_t$, we have
\begin{equation}
P^G_t\left( C < \frac{R_t-\beta^{1/3}t}{\frac{1}{\sqrt{3}} \sqrt{t}} \right) \sim \frac{1}{\sqrt{2\pi}}\int_{C}^{\infty}  e^{  -\frac{1}{2}y^2 +\sum_{k=3}^{\infty} \frac{(-1)^{k+1}}{\beta^{(k+1)/3}} \left(\frac{y}{\sqrt{3t}}\right)^k t} dy\rightarrow 1-\Phi(C).
\end{equation}
\qed\\

\noindent
$\mathbf{Proof \;\;of \;\;Theorem\;\; 1.5}$\\ 
From \cite{Fe51}, we can compute the joint distribution of $B_t$ and $R_t$ under the condition $B_t>0$.
For $0<x<r$,
\begin{equation*}
P(B_t\in dx,R_t\in dr, B_t >0)= \frac{r-x}{t\sqrt{ t}}\cdot \left\{\sum_{k=-\infty}^{\infty} 4k^2\left[ -1 + \left(\frac{2kr-x}{\sqrt{ t}}\right)^2
\right]\phi\left(\frac{2kr-x}{\sqrt{ t}}\right)\right\}
\end{equation*}
\begin{equation}\label{BtRt}
+\sum_{k=1}^{\infty}\left\{ 4k(k-1)\left(\frac{2kr-x}{ t\sqrt{ t}}\right)\phi\left(\frac{2kr-x}{\sqrt{ t}}\right) -4k(k+1) \left(\frac{2kr+x}{t\sqrt{ t}}\right)\phi\left(\frac{2kr+x}{\sqrt{ t}}\right)
\right\}dxdr.
\end{equation}

In the proof of Theorem 1.4, we know that $P^G_t(|R_t/t -\beta^{1/3}|>\epsilon)$ goes to $0$ with an exponential decay. Thus,
$$P_t^G(|B_t/t -\beta^{1/3}|<\epsilon| B_t >0)= P_t^G(|B_t/t -\beta^{1/3}| < \epsilon,|R_t/t -\beta^{1/3}|\leq \epsilon| B_t >0) +o(1)$$
In the  event at the right-hand side, $R_t$ and $B_t$ are both of order $t$, so we only need to consider $k=1$ of the first term in \eqref{BtRt}. We then compute
\[\begin{array}{rl}
 & E\left(1_{\{|\frac{B_t}{t} -\beta^{1/3}| < \epsilon,|\frac{R_t}{t} -\beta^{1/3}|\leq \epsilon\}} e^{-\beta G_t}| B_t >0\right)\\
= & \bigintss_{(\beta^{1/3}-\epsilon)t}^{(\beta^{1/3}+\epsilon)t}  dr \; e^{-\beta\frac{t^2}{r}} \;\bigintss_{(\beta^{1/3}-\epsilon)t}^{r} dx\;\left(\frac{r-x}{t\sqrt{2\pi t}}\right) 4\left[ -1 + \left(\frac{2r-x}{\sqrt{ t}}\right)^2
\right] e^{-\frac{(2r-x)^2}{2t}}\\
= & 4\bigintss_{(\beta^{1/3}-\epsilon)t}^{(\beta^{1/3}+\epsilon)t}  dr \; e^{-\beta\frac{t^2}{r}} \;\bigintss_{(\beta^{1/3}-\epsilon)t}^{r} \;\left(\frac{r-x}{t\sqrt{2\pi t}}\right)  d\left[ (2r-x)e^{-\frac{(2r-x)^2}{2t}}
\right] \\
= & 4\bigintss_{(\beta^{1/3}-\epsilon)t}^{(\beta^{1/3}+\epsilon)t}  dr \; e^{-\beta\frac{t^2}{r}} 
\;\left[ -\frac{(r-x)(2r-x)}{t\sqrt{2\pi t}} e^{-\frac{(2r-x)^2}{2t}}|_{x=(\beta^{1/3}-\epsilon)t}
+ \bigintss_{(\beta^{1/3}-\epsilon)t}^{r} \left(\frac{2r-x}{t\sqrt{2\pi t}}\right)e^{-\frac{(2r-x)^2}{2t}} dx
\right]\\
= & 4\bigintss_{(\beta^{1/3}-\epsilon)t}^{(\beta^{1/3}+\epsilon)t}  dr \; e^{-\beta\frac{t^2}{r}} 
  \;\left[ -\frac{(r-x)(2r-x)}{t\sqrt{2\pi t}} e^{-\frac{(2r-x)^2}{2t}}|_{x=(\beta^{1/3}-\epsilon)t}
  + \frac{1}{\sqrt{2\pi t}}e^{-\frac{r^2}{2t}} - \frac{1}{\sqrt{2\pi t}}e^{-\frac{(2r-x)^2}{2t}}|_{x=(\beta^{1/3}-\epsilon)t}
  \right].
\end{array}\]

Denote $\tilde{x}:=\beta^{1/3}-\epsilon$, $\alpha(r):=-\left(\frac{\beta}{r} +\frac{(2r-\tilde{x})^2}{2}\right)$ is strictly concave. Since $\alpha'(\tilde{x})= \frac{\beta -2\tilde{x}^3}{\tilde{x}^2}<0$, $\alpha$ has the {\it maximizer} at $r=\tilde{x}$ for $\tilde{x}\leq r\leq \tilde{x}+2\epsilon$. The exponent for the first and third term is at most $\alpha(\tilde{x})= -\left(\frac{\beta}{\tilde{x}} +\frac{\tilde{x}^2}{2}\right)$. However, we know the maximizer for $-\left(\frac{\beta}{r} +\frac{r^2}{2}\right)$ is $\beta^{1/3}$, not $\tilde{x}$. 
By applying the Laplace's method again, 
$$E\left(1_{\{|\frac{B_t}{t} -\beta^{1/3}| < \epsilon,|\frac{R_t}{t} -\beta^{1/3}|\leq \epsilon\}} e^{-\beta G_t}| B_t >0\right)$$ 
$$\sim 4\int_{(\beta^{1/3}-\epsilon)t}^{(\beta^{1/3}+\epsilon)t}  dr \; e^{-\beta\frac{t^2}{r}} 
  \;\left[ 
   \frac{1}{\sqrt{2\pi t}}e^{-\frac{r^2}{2t}}
  \right]
  \sim Z_t^G(\{B_t>0\}).$$
We have that $P_t^G(|B_t/t -\beta^{1/3}|<\epsilon| B_t >0)$ goes to 1.

The second part of the proof is devoted for the large deviation principle of the velocity of the polymer end under $\{ B_t>0\}$. It suffices to compute the rate function.\\
Since $R_t$ and $B_t$ are both of order $t$, we only need to consider $k=1$ of the first term in \eqref{BtRt}. We now maximize $-\left(\frac{\beta}{r} +\frac{(2r-\theta)^2}{2}\right)$. It can be found that if $\sqrt[3]{\frac{\beta}{2}}\leq\theta $, the {\it maximizer} is $\theta$ and
$$-J^{\beta}(\theta)=-\frac{\beta}{\theta} -\frac{\theta^2}{2}-g^{**}(\beta).$$
If $\theta \leq \sqrt[3]{\frac{\beta}{2}}$, the {\it maximizer} is $\bar{r}$ and
$$-J^{\beta}(\theta) = -\frac{\beta}{\bar{r}}-\frac{(2\bar{r}-\theta)^2}{2}-g^{**}(\beta),$$
where $\bar{r}$ is the solution of $\beta =2r^2(2r-\theta)$. \qed\\

\noindent
$\mathbf{Proof \;\;of \;\;Theorem\;\; 1.6}$\\ 
We first know
$$P_t^G\left( C< \frac{B_t -\beta^{1/3}t}{\frac{1}{\sqrt{3}}\sqrt{t}} | B_t >0 \right) = P_t^G\left( C< \frac{B_t -\beta^{1/3}t}{\frac{1}{\sqrt{3}}\sqrt{t}} , |R_t/t -\beta^{1/3}|<\epsilon| B_t >0 \right) +o(1).$$
\noindent
We then compute
\[\begin{array}{rl}
 & E(\{ C< \frac{B_t -\beta^{1/3}t}{\frac{1}{\sqrt{3}}\sqrt{t}} , |R_t/t -\beta^{1/3}|<\epsilon\} e^{-\beta G_t}| B_t >0)\\
= & \bigintss_{\beta^{1/3}t+\frac{C\sqrt{t}}{\sqrt{3}}}^{(\beta^{1/3}+\epsilon)t}  dr \; e^{-\beta\frac{t^2}{r}} \;\bigintss_{\beta^{1/3}t+\frac{C\sqrt{t}}{\sqrt{3}}}^{r} dx\;\left(\frac{r-x}{t\sqrt{2\pi t}}\right) 4\left[ -1 + \left(\frac{2r-x}{\sqrt{ t}}\right)^2
\right] e^{-\frac{(2r-x)^2}{2t}}\\
= & 4\bigintss_{\beta^{1/3}t+\frac{C\sqrt{t}}{\sqrt{3}}}^{(\beta^{1/3}+\epsilon)t}  dr \; e^{-\beta\frac{t^2}{r}} 
  \;\left[ \frac{1}{\sqrt{2\pi t}}e^{-\frac{r^2}{2t}}
  -\left(\frac{(r-x)(2r-x)}{t\sqrt{2\pi t}} e^{-\frac{(2r-x)^2}{2t}} + \frac{1}{\sqrt{2\pi t}}e^{-\frac{(2r-x)^2}{2t}}\right)|_{x=\beta^{1/3}t+\frac{C\sqrt{t}}{\sqrt{3}}}
  \right]\\

= & 4\bigintss_{\beta^{1/3}+\frac{C}{\sqrt{3t}}}^{\beta^{1/3}+\epsilon} \;t dr \; e^{-\frac{\beta}{r}t} 
  \;\left[ \frac{1}{\sqrt{2\pi t}}e^{-\frac{r^2}{2}t}
  -\left(\frac{(r-x)t(2r-x)t}{t\sqrt{2\pi t}} e^{-\frac{(2r-x)^2}{2}t} + \frac{1}{\sqrt{2\pi t}}e^{-\frac{(2r-x)^2}{2}t}\right)|_{x=\beta^{1/3}+\frac{C}{\sqrt{3t}}}
  \right].
\end{array}\]
For the second and third term in the last line,
\[\begin{array}{rl}
& -4\bigintss_{\beta^{1/3}+\frac{C}{\sqrt{3t}}}^{\beta^{1/3}+\epsilon} \;t dr \; e^{-\frac{\beta}{r}t} 
  \;
  \left(\frac{(r-x)t(2r-x)t}{t\sqrt{2\pi t}} e^{-\frac{(2r-x)^2}{2}t} + \frac{1}{\sqrt{2\pi t}}e^{-\frac{(2r-x)^2}{2}t}\right)|_{x=\beta^{1/3}+\frac{C}{\sqrt{3t}}}\\
= & e^{-\frac{3}{2}\beta^{3/2}t} O(\frac{1}{\sqrt{t}}).
\end{array}\]
We finally have that
\begin{equation}
P_t^G\left( C< \frac{B_t -\beta^{1/3}t}{\frac{1}{\sqrt{3}}\sqrt{t}} , |R_t/t -\beta^{1/3}|<\epsilon| B_t >0 \right) \rightarrow 1-\Phi(C).
\end{equation}
\qed

\subsection{The discrete model}
The rate function for the large deviation principle for the range of SRW was discussed.

\begin{theorem} (Hamana and Kesten \cite{HK01,HK02}) 
\begin{equation}
I_d(x)=\lim_{n\rightarrow \infty}\frac{-1}{n}\log P\{R_n\geq xn\}
\end{equation}
exists in $[0,\infty]$ for all $x$. $I_d(x)$ is continuous on $[0,1]$ and strictly increasing on $[\gamma_d,1]$, and for $d\geq 2$, $I(x)$ is convex on $[0,1]$. Furthermore,
\begin{equation}
  \begin{aligned}
I_d(x)=0 & \;\;\it{for}\;\; x\leq \gamma_d,\\
0<I_d(x)<\infty & \;\;\it{for}\;\; \gamma_d < x\leq 1,\\
I_d(x)=\infty & \;\;\it{for}\;\; x>1.
   \end{aligned}
\end{equation}
\end{theorem}

Note that $I_d(1)=\log 2d$.
When $d=1$ and $S$ is the SRW, $I(x)=I_1(x)$ can be found explicitly. For $0\leq x\leq 1$,
\begin{equation}
I(x)=\frac{1}{2}(1+x)\log(1+x)+\frac{1}{2}(1-x)\log(1-x).
\end{equation}
The large deviation results were also provided.

\noindent
$\mathbf{Proof \;\;of \;\;Theorem\;\; 1.1}$\\
For any dimension, we set $\tilde{c}(\beta)=\tilde{c}_d(\beta) =\frac{\beta}{\beta+ \log 2d}$. By taking one self-avoiding path, 
\begin{equation}
E\left(e^{-\beta G_n} \mathbf{1}_{R_n\geq \tilde{c}(\beta) n}\right)\geq e^{-\beta n} (2d)^{-n}.
\end{equation}
On the other hand,
\begin{equation}
E\left(e^{-\beta G_n} \mathbf{1}_{R_n < \tilde{c}(\beta) n}\right)\leq e^{-\beta \frac{n}{\tilde{c}(\beta)}}= e^{-(\beta+\log 2d )n}.
\end{equation}
Apply Theorem 2.3 and the Varadhan's lemma we have (i) and the first part of Corollary 1.7.

For $d=1$, we have the explicit rate function $I$. Denote $\psi(c) =\beta/c+ I(c)$. Note that $\psi(\tilde{c})=\beta+ \log 2 + I(\tilde{c})\geq \psi(1)$. Since $\psi$ is differentiable and strictly convex, the minimizer $c^*$ is between $\tilde{c}$ and $1$ by the mean value theorem. The infimum is obtained at $c^*(\beta)$, which is the solution of 
\begin{equation}
\frac{\beta}{c^2} =I '(c) =\frac{1}{2} \log \frac{1+c}{1-c} .
\end{equation}
The monotonicity can be derived since
\begin{equation}
\frac{dc^* (\beta)}{d\beta}= \frac{\sigma^*(\beta)^2}{c^*(\beta)^2}>0
\end{equation}
Recall that $\frac{1}{ \sigma^{*2} (\beta)} = \left( \frac{\beta}{\theta}+I(\theta)\right)''|_{\theta=c^*(\beta)}= \frac{2\beta}{c^{*3}(\beta)} + \frac{1}{1-c^*(\beta)^2}$.
\qed\\

\noindent
$\mathbf{Proof \;\;of \;\;Theorem\;\; 1.2}$\\
We follow the idea in Zoladek \cite{Zo87}. From R{\'e}v{\'e}sz\cite{Re05}, Chapter 2, denote $-M_n^{-}= \min_{0\leq i\leq n} S_i$ and $M_n^{+}= \min_{0\leq i\leq n} S_i$, for $L \leq 0\leq U, L<U, L\leq X\leq U$,
\begin{equation*}
\begin{array}{l}
j(L,U,X):=P(L< -M_n^{-} \leq M_n^+ < U , S_n =X)\\
=\sum_{k=-\infty }^{\infty} P(S_n =X+ 2k(U-L))- \sum_{k=-\infty }^{\infty} P(S_n =2U-X+ 2k(U-L))
\end{array}
\end{equation*}
We then have
\begin{equation*}
\begin{array}{l}
P(L= -M_n^{-} \leq M_n^+ = U , S_n =X)\\
=(j(L-1,U+1,X)-j(L,U+1,X)) -(j(L-1,U,X)-j(L,U,X))
\end{array}
\end{equation*}

From previous discussions, we know that $R_n$ and $S_n$ are of order $n$, and we only consider the case $S_n>0$. Thus, $P(S_n =-X+ 2(U-L))$ dominates other terms. By Sterling's formula, we have
\begin{equation*}
P(S_n=an)=C^{n}_{(1+a)n/2} \;2^{-n} = \frac{1}{\sqrt{2\pi n a(1-a)}} e^{-nI(2a-1)}\left( 1+O\left(\frac{1}{n}\right)\right)
\end{equation*}
Thus, 
\begin{equation}\label{RnSndensity}
P(R_n =rn , S_n =xn)=  \exp(-nI(2r-x)+const\cdot\log n+O(1) ),
\end{equation}
where $U=un, L=\ell n$, and $r=u-\ell$. We then follow the procedure in the proof of Theorem 1.5. \qed\\

\noindent
$\mathbf{Proof \;\;of \;\;Theorem\;\; 1.3}$\\
We prove the CLT under $S_n >0$ and take $c=c^*(\beta)$ for short. First, from Taylor's expasion,
$\frac{1}{r} = \frac{1}{c} -\frac{1}{c^2}(r-c) +\frac{2}{c^3 \cdot 2}(r-c)^2 +\sum_{k=3}^{\infty}\frac{(-1)^{k}}{c^k}(r-c)^k$. Denote $P^c$ be the probability measure for the simple random walk with drift $c$, namely, $P^c_n(X_1=1)=\frac{1+c}{2}$ and $P^c(X_1=-1)=\frac{1-c}{2}$, we then compute
\[\begin{array}{rcl}
& &E\left( \exp\left(-\beta\frac{n^2}{R_n}\right) \right) \\

& =&\sum_{S}\limits \exp\left(-\beta\frac{n^2}{R_n}\right) \left(\frac{1}{1+c}\right)^{\frac{n+S_n}{2}} \left(\frac{1}{1-c}\right)^{\frac{n-S_n}{2}} \left(\frac{1+c}{2}\right)^{\frac{n+S_n}{2}} \left(\frac{1-c}{2}\right)^{\frac{n-S_n}{2}}\\

& =& E^{c} \left( \exp\left(-\beta n \frac{1}{R_n/n} -\frac{S_n}{2} \log \frac{1+c}{1-c}-\frac{n}{2}\log(1-c^2) \right)\right)\\

& \sim& E^{c} \left( \exp\left(-\beta n [\frac{1}{c} -\frac{1}{c^2}(\frac{R_n}{n}-c) +\frac{2}{c^3 \cdot 2}(\frac{R_n}{n}-c)^2 ] -\frac{\beta }{c^2}S_n -\frac{n}{2}\log(1-c^2) \right)\right) \\

& =& E^{c} \left( \exp\left(\frac{\beta}{c^2}  (R_n-S_n) -\frac{2\beta}{c^3 \cdot 2}n(\frac{R_n}{n}-c)^2  \right)\right) \exp\left(-\frac{2\beta}{c} n-\frac{n}{2}\log(1-c^2) \right)\\

& =& 
E^{c} \left( \exp\left(\frac{\beta}{c^2}  (R_n-S_n) -\frac{2\beta}{c^3 \cdot 2}(\frac{R_n-cn}{\sqrt{n}})^2  \right)\right) \exp\left(ng^*(\beta) \right)
\end{array}
\]
\noindent
Recall that $\frac{1}{2} \log \frac{1+c}{1-c}=\frac{\beta}{c^2}$. In the exponent, the order higher than 2 can be ignored. Under $P^c$, $R_n$ and $S_n$ are both asymptotic to a variable with mean $cn$ and variance $(1-c^2)n$. As in the proof for the continuous counterpart, the polymer measure pushes $R_n$ and $S_n$ together, therefore, 

$$E^{c} \left( \exp\left(\frac{\beta}{c^2}  (R_n-S_n)-\frac{2\beta}{c^3 \cdot 2}\left(\frac{R_n-cn}{\sqrt{n}}\right)^2  \right)\right) $$
$$=  E^{c} \left( \exp\left(-\frac{2\beta}{c^3 \cdot 2}\left(\frac{S_n-cn}{\sqrt{n}}\right)^2   + const\cdot \log n +O(1)\right)\right) $$
$$= \int_{-\infty}^{\infty} \exp \left(-\frac{2\beta}{c^3 \cdot 2}z^2 \right) \frac{1}{\sqrt{2\pi  (1-c^2)}}\exp\left( -\frac{z^2}{2(1-c^2)}\right) dz \;\exp\left( const\cdot \log n +O(1)\right)$$
The rest of proof completes by standard arguments. This shows that $\frac{S_n-c^*(\beta)n}{\sqrt{n}} \Longrightarrow N(0, \sigma^*(\beta)^2)$ under the polymer measure $E_n^G$. \qed

\end{document}